\documentclass[12pt]{article}

\usepackage{amsmath,amsthm,amsfonts,amssymb,amscd}
\usepackage{epsfig}

\newtheorem {Lemma}{Lemma}[section]
\newtheorem {Theorem}{Theorem}[section]

\newenvironment {Proof} {\noindent {\bf Proof.}}{\quad $\blacksquare$\par\vspace{3mm}}

= 0.0in \topmargin = 0.3in \oddsidemargin=0.1in
\def\Dj{\hbox{D\kern-.73em\raise.30ex\hbox{-}
\raise-.30ex\hbox{}}}
\def\dj{\hbox{d\kern-.33em\raise.80ex\hbox{-}
\raise-.80ex\hbox{\kern-.40em}}}

\allowdisplaybreaks [4]

\begin{document}

\title{\bf The Kirchhoff indices and the matching numbers of unicyclic graphs}

\author{Xuli Qi$^{a}$, Bo Zhou$^{b,}$\footnote{Corresponding author. Email: {\tt
zhoubo@scnu.edu.cn}}, Zhibin Du$^c$\\
$^a$Hebei Key Laboratory of Computational Mathematics \\and Applications, and\\
College of Mathematics and Information Science, \\ Hebei Normal University, 
 \\
 Shijiazhuang 050024, P. R. China\\
$^b$Department of Mathematics, South China Normal University, \\
Guangzhou 510631, P. R. China\\
$^c$Department of Mathematics, Tongji University,\\
 Shanghai 200092, P. R. China
}

\date{}

\maketitle

\noindent {\bf Abstract } The Kirchhoff index of a connected graph
is the sum of resistance distances between all unordered pairs of
vertices in the graph. It found considerable applications in a variety of fields.
In this paper, we determine the minimum Kirchhoff index among
the unicyclic graphs
with fixed number of vertices and matching number, and characterize the extremal graphs.

\section{Introduction}

The resistance distance was introduced by Klein and Randi\'{c}
\cite{KR} as a distance function on a graph.
Let $G$ be a simple connected graph with vertex set $V(G)$ and edge set
$E(G)$.
The resistance distance between vertices $u$ and $v$ of $G$, denoted by  $r_G(u,v)$, is defined as
the effective resistance
between nodes $u$ and $v$ of the electrical network for
which nodes correspond to the vertices of $G$ and
each edge of $G$ is replaced by a resistor of unit
resistance (one ohm).

The Kirchhoff index of a connected graph $G$ is defined as
\cite{BBLK}
\[
Kf(G)=\sum_{\{u,v\}\subseteq V(G)}r_G(u,v).
\]
It is also named as total effective resistance \cite{GBS}. This
graph invariant found applications in chemistry, electrical network,
Markov chains, averaging networks, experiment design, and Euclidean
distance embeddings, see \cite{Kl,BBLK,GBS}.

The (ordinary) distance between vertices $u$ and $v$ of a graph $G$, denoted by $d_G(u,v)$, is the length of a shortest path connecting them in $G$.
Recall that the Wiener index of $G$ is defined as \cite{DEG,DGKZ} $W(G)=\sum_{\{u,v\}\subseteq V(G)}d_G(u,v)$.
It has been shown \cite{KR} that $r_G(u,v)\le d_G(u,v)$ with
equality if and only if there is a unique path connecting $u$ and
$v$ in $G$. As a consequence, the Kirchhoff index for a tree is equal to its
Wiener index, which has been extensively studied  (see \cite{DEG}). Thus  the
Kirchhoff index is primarily of interest in the case of
cycle-containing graphs.

Zhou and Trinajsti\'{c} \cite{ZT1,ZT2}
established various lower and upper bounds for the Kirchhoff index,
see also \cite{ZT3}. Among the $n$-vertex connected graphs, Lukovits
et al. \cite{LNT} showed that the complete graph $K_n$ is the
unique graph with minimum Kirchhoff index, and Palacios \cite{Pal} showed
 that the path $P_n$ is the unique graph with maximum Kirchhoff
index. The maximum and minimum Kirchhoff indices among the unicyclic
graphs have been determined by Yang and Jiang  \cite{YJ}, see also \cite{ZD}. 

A matching $M$ of the graph $G$ is a subset of $E(G)$ such that no
two edges in $M$ share a common vertex. A matching $M$ of $G$ is
said to be maximum, if for any other matching $M'$ of $G$, $|M'|\le
|M|$. The matching number of $G$ is the number of edges of a maximum
matching in $G$. For a matching $M$ of a graph $G$, if the vertex
$v\in V(G)$ is incident with an edge of $M$, then $v$ is said to be
$M$-saturated. Moreover, if every vertex of $G$ is $M$-saturated, then $M$
is a perfect matching of $G$.

Zhou and Trinajsti\'{c} \cite{ZT4} determined the graphs with
minimum Wiener index and Kirchhoff index respectively among the
connected graphs with fixed number of vertices and matching number.
Du and Zhou \cite{DZ} determined the graphs with  minimum Wiener
index among the trees and unicyclic graphs respectively with fixed
number of vertices and matching number.

In this paper, we determine the minimum Kirchhoff index among
the unicyclic graphs
with fixed number of vertices and matching number, and characterize the extremal graphs.
It is of interest to point out that among the unicyclic graphs with fixed number of vertices and matching number, the graphs with  minimum Kirchhoff index
are different from those with  minimum Wiener index (see \cite{DZ}). 

\section{Preliminaries and Lemmas}

For a graph $G$ with $v\in V(G)$, $G-v$ denotes the graph resulting
from $G$ by deleting $v$  (and its incident edges).  For an edge $uv$ of the graph $G$ (the complement of $G$,
respectively), $G-uv$ ($G+uv$, respectively) denotes the graph
resulting from $G$ by deleting (adding, respectively)  $uv$.


For $u\in V(G)$, let $Kf_G(u)=\sum\limits_{v\in V(G)}r_G(u,v)$. Then
\[
Kf(G)=\frac{1}{2}\sum_{u\in V(G)}Kf_G(u).
\]

Let $C_n$ be the cycle on $n\ge 3$ vertices, whose vertices are
labeled consecutively by $v_1, v_2, \dots,v_n$.

For two vertices $v_i,v_j \in V(C_n)$ with $i < j$, by Ohm's law, we have

\begin{eqnarray} \label{cycle}
r_{C_n} (v_i, v_j) = \frac{(j - i) \cdot [ n - ( j - i)]} {n}.
\end{eqnarray}
Furthermore, for fixed $n$, $r_{C_n} (v_i, v_j)$ is increasing for $j - i \le \lfloor \frac{n}{2} \rfloor$. For $v_1 \in V (C_n)$, by Eq. (\ref{cycle}), we have
\begin{eqnarray} \label{cycle-vertex}
Kf_{C_n} (v_1) = \sum_{i = 2}^n r_{C_n} (v_1,v_i) = \sum_{i = 2}^n \frac{(i - 1) \cdot [ n - ( i - 1)]} {n} = \frac {n^2-1}{6},
\end{eqnarray}
and thus
\begin{eqnarray} \label{cycle-all}
Kf (C_n) = \frac{1} {2} \cdot n \cdot Kf_{C_n} (v_1) = \frac {n^3 - n}{12}.
\end{eqnarray}

For a unicyclic graph $G$ with the unique cycle $C_k$,  $G-E(C_k)$ consists of
$k$ vertex-disjoint trees $T_1, T_2, \dots, T_{k}$, where $v_i\in V(T_i)$ for
$i=1,2,\dots, k$. These trees are called the branches of $G$, and
$v_i$ is called the root of the branch $T_i$ in $G$ for $i=1, 2, \dots,
k$.

\begin{figure}[h]
  \centering
  \includegraphics[height=4.5cm]{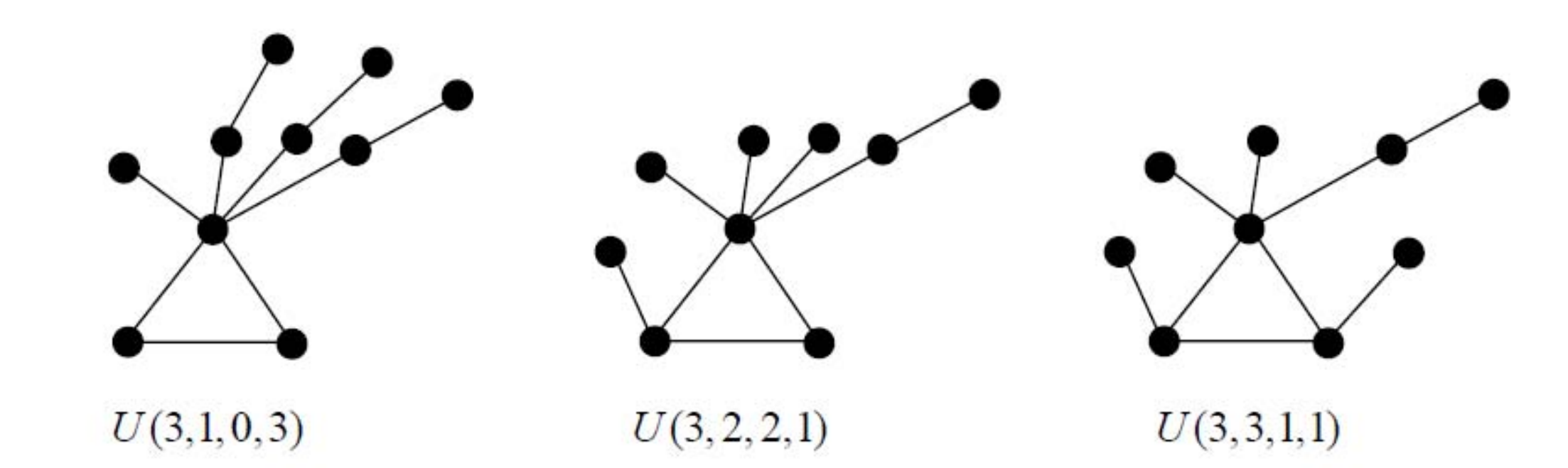}
  \caption{The graphs $U (3,1,0,3)$, $U (3,2,2,1)$ and $U (3,3,1,1)$.} \label{fig1}
\end{figure}

Now we define the graph $U(k,t,i,j)$ which will be used frequently later. For integers $k,t,i,j$ with $k \ge 3$, $k \ge t \ge 0$, $i \ge 0$, $j \ge 0$, let
$U(k,t,i,j)$ be the graph obtained from the cycle $C_k$ as follows:

\begin{enumerate}
  \item[($a$)]
  choose $t$ consecutive vertices in the cycle $C_k$$;$
  \item[($b$)]
  attach $t$ pendent vertices each to one of the $t$ chosen
vertices in (a)$;$
  \item[($c$)]
  attach $i$ pendent vertices and $j$ paths on two vertices to
a central vertex of the $t$ chosen vertices in (a).
\end{enumerate}
Clearly, $U(k,t,i,j)$ has $k+t+i+2j$ vertices. 
%
%
In particular, let $U (k,t) = U (k,t,0,0)$ for integers $k,t$ with $k \ge 3$ and
$k \ge t \ge 0$.
For example, $U (3,1,0,3)$, $U (3,2,2,1)$ and $U (3,3,1,1)$ are shown in Fig. \ref{fig1}.

Let $d_G (v)$ be the degree of $v$ in $G$.

For integers $n$ and $m$ with $2\le m\le\lfloor\frac{n}{2}\rfloor$,
let $\mathbb{U}(n,m)$ be the set of unicyclic graphs with $n$
vertices and matching number $m$.
For integer $m \ge 2$, we can partition $\mathbb{U}(2m,m)\setminus \{C_{2m}\}$ into two subsets as follows:
\begin{enumerate}
\item[$(i)$]
the set of
graphs of maximum degree three in $\mathbb{U}(2m,m)$ obtainable by
attaching some pendent vertices to a cycle, which is denoted by $\mathbb{U}_1(m)$$;$
\item[$(ii)$]
the set of
graphs in $\mathbb{U}(2m,m)$ containing some pendent vertex whose
unique neighbor is of degree two, which is denoted by $\mathbb{U}_2(m)$$.$
\end{enumerate}

\subsection{The Kirchhoff index of graphs in $\mathbb{U}_1(m)$ with small $m$}

First we want to determine the minimum
Kirchhoff index among the graphs in $\mathbb{U}_1(m)$ with $2 \le m \le 8$.


\begin{Lemma}  \label{U-2-m-a}
Let $G\in \mathbb{U}_1(m)$  with the unique cycle $C_k$ and $t$ pendent
vertices, where $k +t = 2m$, $k\ge 3$ and $k \ge t \ge 1$.
\begin{enumerate}
\item[$(i)$]
For $t=1,2,3,k-4,k-2,k$,
$(k,t)=(10,4)$, $(k,t)=(11,5)$, or $(k,t)=(12,4)$, we have
\[
Kf(G)\ge
\frac{1}{12}\left(k^3+2k^2t+12kt-k+2t^3+12t^2-16t+\frac{t^2-t^4}{k}\right)
 \]
with equality if and only if $G \cong U (k,t)$.
\item[$(ii)$]
For integers $k,t$ with $k \ge 3$,
$k \ge t \ge 1$, and $v \in V(G)$, we have
\[
Kf_G (v) \ge f(k,t)
\]
with equality if and only if $G \cong U(k,t)$, and $v$ is a central
vertex of the $t$ vertices of degree three in $U(k,t)$, where
\[
f(k,t)=
\begin{cases}
\frac{1}{12}(2k^2+3t^2+12t-5-\frac{t^3-t}{k})&\text{if $t$ is odd},\\
\frac{1}{12}(2k^2+3t^2+12t-2-\frac{t^3+2t}{k})&\text{if $t$ is even}.
\end{cases}
\]
\end{enumerate}
\end{Lemma}

\begin{Proof}
First we prove (i). The cases $t=1,k-2,k$ are trivial. Suppose that $t\not=1,k-2,k$.

Let $S(G)=\{v\in V(C_k):d_G(v)=3\}$, and let $\sigma
(G)=\sum\limits_{\{v_i,v_j\}\subseteq S(G)} r_G(v_i,v_j)$. Clearly, $| S (G) | = t$.

If $t=2$, say $S(G)=\{v_1,v_s\}$, then
\[
\sigma (G)=r_G(v_1,v_s)\ge \frac{1\cdot (k-1)}{k} = \sigma(U(k,2))
\]
with equality if and only if $v_1$ and $v_s$ are adjacent in $G$, i.e., $G
\cong U(k,2)$.

If $t\ge 3$, then by Eq. (\ref{cycle}), we have
\begin{eqnarray} \label{U(k,t)}
\sigma (U(k,t))
&=& r_G (v_1,v_2) + r_G (v_1,v_3) + \cdots + r_G (v_1,v_t) \nonumber \\
&& + r_G (v_2,v_3) + r_G (v_2,v_4) + \cdots + r_G (v_2,v_t) \nonumber \\
&& + \cdots + r_G (v_{t - 1},v_t) \nonumber \\
&=& \sum_{i = 1}^{t - 1} \sum_{j = i + 1}^{t} r_G (v_i,v_j) \nonumber \\
&=& \sum_{i = 1}^{t - 1} \sum_{j = i + 1}^{t} \frac{( j - i ) \cdot [ k - (j - i)]}{k} \nonumber \\
&=& \sum_{i = 1}^{t - 1} \sum_{j = 1}^{t - i} \frac{j \cdot (k-j)}{k} \nonumber \\
&=& \sum_{i = 1}^{t - 1} \sum_{j = 1}^{i} \frac{j \cdot (k-j)}{k} \nonumber \\
&=&=\frac{1}{12k}t(t-1)(t+1)(2k-t).
\end{eqnarray}

Suppose that $t=3$. Then $k$ is odd as $G$ has perfect matching. By
symmetry, we may assume that $S(G)=\{v_1,v_i,v_j\}$ with $1 < i <
j$, and $d_{G}(v_1,v_i)\le d_{G}(v_1,v_j)$. Obviously, $i\le
\frac{k+1}{2}$. If $j\le \frac{k+1}{2}$, then note that $d_G (v_1,v_j) \ge 2$, and by Eq. (\ref{cycle}), we have
\begin{eqnarray*}
\sigma (G)&=&r_G(v_1,v_i) + r_G(v_1,v_j) + r_G(v_i,v_j)\\
&\ge& \frac{1\cdot (k-1)}{k} + \frac{2\cdot (k-2)}{k} + \frac{1\cdot
(k-1)}{k} \\
&=& \sigma(U(k,3))
\end{eqnarray*}
with equality if and only if $i=2$ and $j=3$, i.e., $G \cong
U(k,3)$. 
If $j=k$, then we
have $i=2$ since $d_{G}(v_1,v_i)\le d_{G}(v_1,v_j)$, i.e., $G
\cong U(k,3)$. Note that $j\ne k-1$ as $G$ has perfect matching. If
$\frac{k+3}{2}\le j\le k-2$, then $d_G (v_1,v_j) \ge 3$, and by Eq. (\ref{cycle}), we have
\begin{eqnarray*}
\sigma (G)&=&  r_G(v_1,v_i) + r_G(v_1,v_j) + r_G(v_i,v_j)\\
&\ge& \frac{1\cdot (k-1)}{k} + \frac{3\cdot (k-3)}{k} + \frac{1\cdot
(k-1)}{k}\\
&>& \frac{1\cdot (k-1)}{k} + \frac{2\cdot (k-2)}{k} + \frac{1\cdot
(k-1)}{k}\\
&=& \sigma(U(k,3)).
\end{eqnarray*}
Now it follows that $\sigma (G)\ge \sigma (U(k,3))$ with equality if and only if
$G \cong U(k,3)$.

Suppose that $t=k-4\ge 4$. Suppose to the contrary that $G \not\cong
U(k,k-4)$. Then there are two pairs of adjacent vertices of degree
two on the cycle $C_k$ in $G$, separated by $a\ge 1$ consecutive
vertices $v_{i_1},v_{i_2},\ldots,v_{i_a}$ of degree three and $b\ge
1$ consecutive vertices $v_{j_1},v_{j_2},\ldots,v_{j_b}$ of degree
three on the cycle $C_k$, where $d_{C_k}(v_{i_1},v_{j_1})=3$,
$d_{C_k}(v_{i_a},v_{j_b})=3$, and $a+b=k-4$. Assume that $a\ge b$.
Denote by $w$ the pendent neighbor of $v_{j_1}$ in $G$. Consider
$G'=G- \{ v_{j_1} w \} + \{ vw \}$, where $v$ is the neighbor of $v_{i_1}$ with
degree two on the cycle. Note that
$S(G)=\{v_{i_1},v_{i_2},\dots,v_{i_a},v_{j_1},v_{j_2},\dots,v_{j_b}\}$
and
$S(G')=\{v_{i_1},v_{i_2},\dots,v_{i_a},v,v_{j_2},\dots,v_{j_b}\}$.
If $b\ge 2$, then by Eq. (\ref{cycle}), we have
\begin{eqnarray*}
\sum_{s = 1}^{a} r_G(v_{j_1},v_{i_s})-\sum_{s = 1}^{a}  r_{G'}(v,v_{i_s})
=\frac{1}{k}\sum_{i=3}^{a+2}i(k-i)-\frac{1}{k}\sum_{i=1}^{a}i(k-i),
\end{eqnarray*}
\begin{eqnarray*}
\sum_{s = 2}^{b}  r_G(v_{j_1},v_{j_s})-\sum_{s = 2}^{b} r_{G'}(v,v_{j_s})
=\frac{1}{k}\sum_{i=1}^{b-1}i(k-i)-\frac{1}{k}\sum_{i=3}^{b+1}i(k-i),
\end{eqnarray*}
and thus
\begin{eqnarray*}
&&\sigma (G)-\sigma (G')\\
&=& \sum_{x \in S (G) \setminus \{v_{j_1}\}} r_G(v_{j_1},x) - \sum_{x \in S (G') \setminus \{v\}} r_{G'}(v,x) \\
&=&\left( \sum_{s = 1}^{a} r_G(v_{j_1},v_{i_s}) +
\sum_{s = 2}^{b} r_G(v_{j_1},v_{j_s}) \right)\\
&&-\left( \sum_{s = 1}^{a} r_{G'}(v,v_{i_s}) +
\sum_{s = 2}^{b} r_{G'}(v,v_{j_s}) \right)\\
&=&\left(\sum_{s = 1}^{a} r_G(v_{j_1},v_{i_s})-\sum_{s = 1}^{a} r_{G'}(v,v_{i_s})\right)\\
&&+ \left(\sum_{s = 2}^{b} r_G(v_{j_1},v_{j_s})-\sum_{s = 2}^{b} r_{G'}(v,v_{j_s})\right)\\
&=& \left(
\frac{1}{k}\sum_{i=3}^{a+2}i(k-i)-\frac{1}{k}\sum_{i=1}^{a}i(k-i)
\right) \\
&&+ \left(
\frac{1}{k}\sum_{i=1}^{b-1}i(k-i)-\frac{1}{k}\sum_{i=3}^{b+1}i(k-i)
\right) \\
&=&\frac{4}{k}(a-b+1)>0.
\end{eqnarray*}
If $b=1$, then by similar arguments as above, we have $\sigma
(G)-\sigma (G')=\frac{4a}{k}>0$. Thus $\sigma (G)>\sigma (G')$ for
$b\ge 1$. By repeating the transformation from $G$ to $G'$, we may finally  get
$\sigma (G) > \sigma (U(k,k-4))$. Thus if $t=k-4$, then  $\sigma
(G)\ge \sigma(U(k,k-4))$ with equality if and only if $G \cong
U(k,k-4)$.

Suppose that $(k,t)=(10,4)$. Then there are exactly four possibilities for $G$, and by suitable labeling, we may assume that
$S(G)=\{v_1,v_2,v_3,v_4\}$, $\{v_1,v_2,v_3,v_6\}$, $\{v_1,v_2,v_5,v_6\}$, or
$\{v_1,v_2,v_5,v_8\}$. By direct calculation, we have
\[
\sigma (G) =
\begin{cases}
8 & \text{if $S (G) = \{v_1,v_2,v_3,v_4\}$},\\
\frac{52}{5} & \text{if $S (G) = \{v_1,v_2,v_3,v_6\}$},\\
\frac{56}{5} & \text{if $S (G) = \{v_1,v_2,v_5,v_6\}$},\\
12 & \text{if $S (G) = \{v_1,v_2,v_5,v_8\}$},
\end{cases}
\]
and thus $\sigma
(G) \ge 8$ with equality if and only if $G \cong U (10,4)$.

Suppose that $(k,t)=(11,5)$. Then there are exactly five possibilities for $G$, and by suitable labeling, we may assume that
$S(G) = \{v_1,v_2,v_3,v_4,v_5\}$, $\{v_1,v_2,v_3,v_4,v_7\}$, $\{v_1,v_2,v_3,v_6,v_7\}$,
$\{v_1,v_2,v_3,v_6,v_9\}$, or $\{v_1,v_2,v_5,v_8,v_9\}$. By direct calculation, we have
\[
\sigma (G) =
\begin{cases}
\frac{170}{11} & \text{if $S (G) = \{v_1,v_2,v_3,v_4,v_5\}$},\\
\frac{202}{11} & \text{if $S (G) = \{v_1,v_2,v_3,v_4,v_7\}$},\\
\frac{218}{11} & \text{if $S (G) = \{v_1,v_2,v_3,v_6,v_7\}$},\\
\frac{226}{11} & \text{if $S (G) = \{v_1,v_2,v_3,v_6,v_9\}$},\\
\frac{234}{11} & \text{if $S (G) = \{v_1,v_2,v_5,v_8,v_9\}$},
\end{cases}
\]
and thus $\sigma (G) \ge \frac{170}{11}$ with equality if and only if $G \cong U (11,5)$.

Suppose that $(k,t)=(12,4)$. Then there are exactly eight possibilities for $G$, and by suitable labeling, we may assume that
$S(G) = \{v_1,v_2,v_3,v_4\}$, $\{v_1,v_2,v_3,v_6\}$, $\{v_1,v_2,v_3,v_8\}$,
$\{v_1,v_2,v_5,v_6\}$, $\{v_1,v_2,v_7,v_8\}$, $\{v_1,v_2,v_5,v_8\}$,
$\{v_1,v_2,v_5,v_{10}\}$, or $\{v_1,v_4,v_7,v_{10}\}$. By direct calculation, we have
\[
\sigma (G) =
\begin{cases}
\frac{25}{3} & \text{if $S (G) = \{v_1,v_2,v_3,v_4\}$},\\
\frac{34}{3} & \text{if $S (G) = \{v_1,v_2,v_3,v_6\}$},\\
\frac{37}{3} & \text{if $S (G) = \{v_1,v_2,v_3,v_8\}$},\\
\frac{37}{3} & \text{if $S (G) = \{v_1,v_2,v_5,v_6\}$},\\
\frac{41}{3} & \text{if $S (G) = \{v_1,v_2,v_7,v_8\}$},\\
14 & \text{if $S (G) = \{v_1,v_2,v_5,v_8\}$},\\
\frac{41}{3} & \text{if $S (G) = \{v_1,v_2,v_5,v_{10}\}$},\\
15 & \text{if $S (G) = \{v_1,v_4,v_7,v_{10}\}$},
\end{cases}
\]
and thus $\sigma (G) \ge \frac{25}{3}$ with equality if and only if $G \cong U (12,4)$.

Combining all the above cases, and by Eq. (\ref{U(k,t)}), we can deduce that
\[
\sigma (G)\ge \sigma (U(k,t)) = \frac{1}{12k}t(t-1)(t+1)(2k-t)
\]
with equality if and only if  $G \cong U(k,t)$ for
$t=1,2,3,k-4,k-2,k$, $(k,t)=(10,4)$, $(k,t)=(11,5)$, or
$(k,t)=(12,4)$. For $1\le i\le k$ with $d_G(v_i)=3$, let $u_i$ be the pendent neighbor of $v_i$ in $G$.

By Eqs. (\ref{cycle-all}) and (\ref{cycle-vertex}), we have
\begin{eqnarray*}
Kf(G)&=&\sum_{\{v_i,v_j\}\subseteq V(C_k)} r_G(v_i,v_j) + \sum_{u_i\in V(G)\setminus V(C_k)}\sum_{v_j\in V(C_k)}r_G(u_i,v_j)\\
     &&+\sum_{\{u_i,u_j\}\subseteq V(G )\setminus V(C_k)} r_G(u_i,u_j)\\
     &=&\frac {k^3-k}{12}+\sum_{v_i \in S(G)}\sum_{v_j\in V(C_k)}(1+r_G(v_i,v_j)) + \sum_{\{v_i,v_j\}\subseteq S(G)}(2+r_G(v_i,v_j))\\
     &=&\frac {k^3-k}{12}+\sum_{v_i\in S(G)}(k + Kf_{C_k}(v_i)) + 2{t\choose 2}+\sum_{\{v_i,v_j\}\subseteq S(G)}r_G(v_i,v_j)\\
     &=&\frac{k^3-k}{12}+t\left(k+\frac{k^2-1}{6}\right)+2{t\choose
     2}+\sigma (G)\\
     &\ge&\frac{k^3-k}{12}+t\left(k+\frac{k^2-1}{6}\right)+2{t\choose
     2}+\frac{1}{12k}t(t-1)(t+1)(2k-t)\\
     &=&\frac{1}{12}\left(k^3+2k^2t+12kt-k+2t^3+12t^2-16t+\frac{t^2-t^4}{k}\right)
\end{eqnarray*}
with equality if and only if $G \cong U(k,t)$ for
$t=1,2,3,k-4,k-2,k$, $(k,t)=(10,4)$, $(k,t)=(11,5)$, or
$(k,t)=(12,4)$.


Next we prove (ii).
Let $v \in V(G)$. For $v_i \in V (C_k)$, clearly $Kf_G (v_i^*) - Kf_G (v_i) = 2m - 2 > 0$, where $d_G (v_i) = 3$, and $v_i^*$ is the unique neighbor of $v_i$ in $G$ outside $C_k$. Thus we may assume that $v = v_i \in V (C_k)$.
By Eq. (\ref{cycle}), it is easily seen that
\begin{eqnarray*}
\sum_{v_j\in S(G)}r_G(v_i,v_j)
&\ge&\begin{cases}0+2\sum\limits_{i=1}^{(t-1)/2}\frac{i\cdot (k-i)}{k}&\text{if $t$ is odd}\\
                         0+2\sum\limits_{i=1}^{(t-2)/2}\frac{i\cdot (k-i)}{k}+\frac{t/2\cdot (k-t/2)}{k}&\text{if $t$ is even}
            \end{cases}\\
            &=&\begin{cases}\frac{1}{12}(3t^2-3-\frac{t^3-t}{k})&\text{if $t$ is odd}\\
                         \frac{1}{12}(3t^2-\frac{t^3+2t}{k})&\text{if $t$ is even}
            \end{cases}
\end{eqnarray*}
with equality if and only if the $t$ vertices in $S(G)$ are
consecutive on $C_k$, i.e., $G \cong U(k,t)$, and $v_i$ is a central
vertex of the $t$ vertices of degree three in $U(k,t)$. For $v_i\in V(C_k)$, by Eq. (\ref{cycle-vertex}), we have
\begin{eqnarray*}
Kf_G(v_i)&=&Kf_{C_k}(v_i)+\sum_{u_j\in V(G)\setminus V(C_k)}r_G(v_i,u_j)\\
         &=&\frac{k^2-1}{6}+\sum_{v_j \in S(G)}(1+r_G(v_i,v_j))\\
         &=&\frac{k^2-1}{6}+t+\sum_{v_j\in S(G)}r_G(v_i,v_j)\\
         &\ge&\begin{cases}\frac{k^2-1}{6}+t+\frac{1}{12}(3t^2-3-\frac{t^3-t}{k})&\text{if $t$ is odd}\\
                         \frac{k^2-1}{6}+t+\frac{1}{12}(3t^2-\frac{t^3+2t}{k})&\text{if $t$ is even}
            \end{cases}\\
         &=&\begin{cases}\frac{1}{12}(2k^2+3t^2+12t-5-\frac{t^3-t}{k})&\text{if $t$ is odd}\\
                         \frac{1}{12}(2k^2+3t^2+12t-2-\frac{t^3+2t}{k})&\text{if $t$ is even}
            \end{cases}\\
&=&f(k,t)
\end{eqnarray*}
with equality if and only if $G \cong U(k,t)$, and $v_i$ is a central
vertex of the $t$ vertices of degree three in $U(k,t)$.
\end{Proof}

If $G\in \mathbb{U}_1(m)$ with the unique cycle $C_k$ and $t$ pendent
vertices, where $2 \le m \le 8$, then $t=1,2,3,k-4,k-2,k$,
$(k,t)=(10,4)$, $(k,t)=(11,5)$, or $(k,t)=(12,4)$. Now by Lemma
\ref{U-2-m-a} (i), we have

\begin{Lemma} \label{U-2-m}
If $G$ is a graph in $\mathbb{U}_1(m)$ with the minimum Kirchhoff index,
where $2 \le m \le 8$, then $G \cong U(k,t)$ with $k +t = 2m$, $k\ge 3$ and $k \ge t \ge 1$.
\end{Lemma}

\subsection{The Kirchhoff index of graphs in $\mathbb{U}_2(m)$ with small $m$}

The following result will be useful for comparing the Kirchhoff
indices of graphs.

For simplicity, let $|G|=|V(G)|$ for a graph $G$.

\begin{Lemma}  \label{formula}\textnormal{\cite{ZD}}
Let $G$ and $H$ be two connected graphs with $u\in V(G)$ and $w\in
V(H)$. Let $GuH$ be the graph obtained from $G$ and $H$ by
identifying $u\in V(G)$ with $w\in V(H)$. Then
\[
Kf(GuH)=Kf(G)+Kf(H)+(|H|-1)\,Kf_{G}(u)+(|G|-1)\,Kf_{H}(w).
\]
\end{Lemma}

Let $P_n$ be the path on $n$ vertices.

If $u$ is a pendent vertex being adjacent to a vertex $v$ of degree
two in the graph $G$, then the path of $G$ induced by the vertices
$u$ and $v$ is said to be a pendent $P_2$ of $G$.
Clearly, every graph in $\mathbb{U}_2 (m)$ has at least one pendent
$P_2$.

For a given graph $G \in \mathbb{U}_2 (m)$, starting from $G$, deleting the pendent
$P_2$'s repeatedly, until there is no pendent
$P_2$, the resulting graph is denoted by $\bar{G}$. Let $\bar{n} = |\bar{G}|$. Clearly, $\bar{G} \in \mathbb{U}_1 (\frac{\bar{n}} {2}) \cup \{C_{\bar{n}}\}$.


Now we determine the minimum
Kirchhoff index among the graphs in $\mathbb{U}_2(m)$ with $3 \le m \le 8$.

\begin{Lemma}  \label{U-1-m}
If $G$ is a graph in $\mathbb{U}_2(m)$ with the minimum Kirchhoff index,
where $3\le m\le8$, then $G \cong U(k,t,0,j)$ with $k+t+2j =2m$, $k \ge
3$, $k \ge t \ge 1$ and $j \ge 1$.
\end{Lemma}

\begin{Proof}
Let $G \in \mathbb{U}_2 (m)$, and $k$ be the length of the unique cycle of $G$.

Denote by the deleting process from $G$ to $\bar{G}$ as follows:
$$
G = G_1 \rightarrow G_2 \rightarrow \dots \rightarrow G_{r-1} \rightarrow G_r = \bar{G},
$$
where $G_{i+1}$ is the (unicyclic) graph obtained from $G_i$ by deleting a pendent $P_2$, where $1 \le i \le r - 1$.
Note that $\bar{n} + 2 (r - 1) = 2 m$.


Recall that $G_r = \bar{G} \in \mathbb{U}_1 (\frac{\bar{n}} {2}) \cup \{C_{\bar{n}}\}$. By Lemma \ref{U-2-m}, we have $Kf (G_{r}) \ge Kf (U (k,t,0,0))$, where $k + t = \bar{n}$. Moreover, by
Lemma \ref{U-2-m-a} (ii) and Lemma \ref{formula}, we have $Kf (G_{r - 1}) \ge Kf (U (k,t,0,1))$ with equality if and only if $G_{r - 1} \cong U (k,t,0,1)$.
Again by Lemma \ref{U-2-m-a} (ii) and Lemma \ref{formula}, we have $Kf (G_{r - 2}) \ge Kf (U (k,t,0,2))$ with equality if and only if $G_{r - 2} \cong U (k,t,0,2)$. Repeating the arguments, finally we can deduce that $Kf (G) = Kf (G_1) \ge Kf (U (k,t,0,r - 1))$ with last equality if and only if $G \cong U (k,t,0,r - 1)$.

Then the result follows easily.
\end{Proof}

\subsection{The effect on the Kirchhoff index of graphs under the deletion of some vertices}

First we introduce a unicyclic graph.

Let $U_{n,m}= U (5,1,n-2m,m-3)$, where $3\le m\le
\lfloor\frac{n}{2}\rfloor$, see Fig \ref{U(n,m)}. It is easily checked that
\begin{eqnarray} \label{eq-Unm}
Kf (U_{n,m}) = n^2+nm-5n-3m+4.
\end{eqnarray}

\begin{figure}[h]
  \centering
  \includegraphics[height=4.5cm]{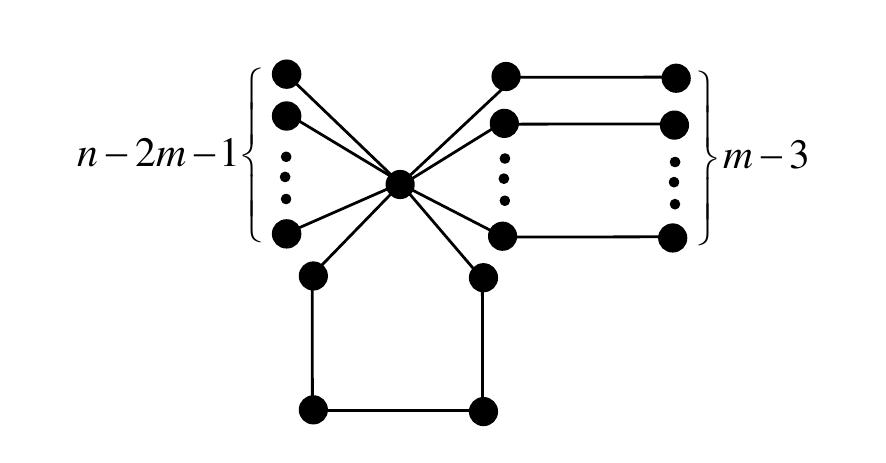}
  \caption{The graph $U_{n,m}$.}\label{U(n,m)}
\end{figure}

Next we establish a lower bound of $Kf_{G}(u)$, where $G\in
\mathbb{U}(n,m)$ and $u \in V(G)$.

\begin{Lemma}\label{degree}
Let $G\in \mathbb{U}(n,m)$ with the unique cycle $C_k$, where $n \ge 6$, $m \ge 3$, $k \ge 3$. If
$T_i \cong P_1$ or $P_2$ for $2\le i \le k$, then for $u\in V(T_1)$,
\[
Kf_{G}(u)\ge n+m-4
\]
with equality if and only if $G \cong U_{n,m}$, and $u$ is the vertex
of maximum degree in $U_{n,m}$.
\end{Lemma}

\begin{Proof}
Let $M$ be a maximum matching of $G$. First we establish an upper bound of $d_G(u)$. Let
\[
A_1 = \{ xy \in E(G) \setminus M : \mbox{either $x = u$ or $y = u$}
\},
\]
\[
A_2 = \{ xy \in E(G) \setminus M : \mbox{$x,y \ne u$ and $xy \in E(C_k)$}  \},
\]
\[
A_3 = \{ xy \in E(G) \setminus M : \mbox{$x,y \ne u$ and $xy \not\in E(C_k)$} \}.
\]
Clearly, $A_1$, $A_2$, $A_3$ are pairwise disjoint, and
$E(G) \setminus M = A_1 \cup A_2 \cup A_3$.
Thus
\begin{equation}\label{outside}
|E(G) \setminus M| = n - m = |A_1| + |A_2| + |A_3|.
\end{equation}
Note that
\begin{enumerate}
\item[(a)]
$|A_1| \ge d_G (u) - 1$ with equality if and only if $u$ is
$M$-saturated$;$
\item[(b)]
$|A_2| \ge \lfloor \frac{k - 2} {2} \rfloor$ if $u$ lies on the unique cycle $C_k$ of $G$, and\\ $|A_2| \ge \lfloor \frac{k + 1} {2} \rfloor$ if
$u$ lies outside the unique cycle $C_k$ of $G$$;$
\item[(c)]
$|A_3| \ge 0$.
\end{enumerate}
It follows from Eq. (\ref{outside}) that if $u$ lies on the unique cycle $C_k$ of $G$, then $n - m \ge
(d_G (u) - 1) + \left \lfloor \frac{k - 2} {2} \right \rfloor$,
i.e.,
\begin{eqnarray} \label{upp-1}
d_G (u) \le n - m + 1 - \left \lfloor \frac{k - 2} {2} \right
\rfloor
\end{eqnarray}
with equality if and only if the corresponding equalities in (a), (b), (c) hold,
while if $u$ lies outside the unique cycle $C_k$ of $G$, then $n - m \ge
(d_G (u) - 1) + \left \lfloor \frac{k + 1} {2} \right \rfloor$,
i.e.,
\begin{eqnarray} \label{upp-2}
d_G (u) \le n - m + 1 - \left \lfloor \frac{k + 1} {2} \right
\rfloor
\end{eqnarray}
with equality if and only if  the corresponding equalities in (a), (b), (c) hold.


%
%

%

\noindent {\bf Case 1.} $u$ lies on the unique cycle $C_k$ of $G$.

\noindent {\bf Subcase 1.1. } $k$ is odd and $T_i \cong P_1$ for $2\le i \le k$.

By Eq. (\ref{cycle-vertex}) and inequality (\ref{upp-1}), we have
\begin{eqnarray*}
Kf_{G}(u)&=&Kf_{C_k}(u)+ \sum_{x \in V (G) \setminus V (C_k)} r_G (u,x)\notag\\
       &\ge&\frac{k^2-1}{6}+[(d_{G}(u)-2)+2(n-k-d_{G}(u)+2)]\notag\\
       &=&-d_{G}(u)+\frac{1}{6}k^2-2k+2n+\frac{11}{6}\notag\\
       &\ge&-\left( n - m + 1 -   \frac{k - 3} {2}
        \right)+\frac{1}{6}k^2-2k+2n+\frac{11}{6}\notag\\
       &=&\frac{1}{6}(k^2-9k+6n+6m-4) 
\end{eqnarray*}
with equality if and only if $G \cong U (k,1,n-2m,m-\frac{k+1}{2})$
with odd $k$, and $u$ is the vertex of maximum degree in
$U(k,1,n-2m,m-\frac{k+1}{2})$.

\noindent {\bf Subcase 1.2.}
$k$ is even and $T_i \cong P_1$ for $2\le i \le k$, or
there is at least one of $T_i$ such that $T_i \cong P_2$ for $2\le i
\le k$.

Obviously, $r_G(u,v)\ge 1+\frac{2\cdot (k-2)}{k}\ge 2$,
where $v$ is the unique pendent neighbor of $v_i$ with $3\le i \le k-1$ if
$T_i \cong P_2$. On the other hand, we also note that if $k$ is odd, then $T_i \cong
P_2$ for some $i$ with $2\le i \le k$, and thus either $|A_2| >
\lfloor \frac{k - 2} {2} \rfloor = \frac{k - 3} {2}$ or $|A_3| > 0$, by Eq. (\ref{outside}), we have $n - m \ge
(d_G (u) - 1) + \frac{k - 3} {2} + 1$,
i.e.,
\begin{eqnarray} \label{upper}
d_G (u) \le
n - m  - \frac{k - 3} {2}.
\end{eqnarray}
Let $a$ be the number of pendent vertices attached to $v_2$ or
$v_k$ in $G$, where $0\le a \le 2$. Then by Eq. (\ref{cycle-vertex}), and inequalities (\ref{upp-1}) and (\ref{upper}), we have
\begin{eqnarray*}
Kf_{G}(u)&=&Kf_{C_k}(u)+ \sum_{x \in V (G) \setminus V (C_k)} r_G (u,x)\\
       &\ge&\frac{k^2-1}{6}+\left[(d_{G}(u)-2)+\left(1+\frac{1\cdot (k-1)}{k}\right)a+2(n-k-a-d_{G}(u)+2)\right]\\
       &=&-\frac{a}{k}-d_{G}(u)+\frac{1}{6}k^2-2k+2n+\frac{11}{6}\\
       &\ge&\begin{cases}-\frac{2}{k}-\left(n-m + 1 -\frac{k-2}{2}\right)+\frac{1}{6}k^2-2k+2n+\frac{11}{6} &\text{if $k$ is even}\\
                         -\frac{2}{k}-\left(n-m -\frac{k-3}{2} \right)+\frac{1}{6}k^2-2k+2n+\frac{11}{6} &\text{if $k$ is odd}
            \end{cases}\\
       &=&\begin{cases}\frac{1}{6}(k^2-9k+6n+6m-1-\frac{12}{k})&\text{if $k$ is even}\\
                        \frac{1}{6}(k^2-9k+6n+6m +2-\frac{12}{k})&\text{if $k$ is odd}
            \end{cases}\\
       &\ge&\frac{1}{6}(k^2-9k+6n+6m-4).
\end{eqnarray*}
If $Kf_{G}(u) = \frac{1}{6}(k^2-9k+6n+6m-4)$, then $a=2$,
$k=4$ and $d_G(u)=n-m + 1 -\frac{k-2}{2}=n-m$. However, $a=2$ and $k=4$
imply that either $|A_2| > 1$ or $|A_3| > 0$, and thus by Eq. (\ref{outside}), we have
$n-m > d_G(u)$.
Therefore
$Kf_{G}(u)>\frac{1}{6}(k^2-9k+6n+6m-4)$.

\noindent {\bf Case 2.} $u$ lies outside the unique cycle $C_k$ of $G$.

Note that for $v\not\in V(C_k)\cup V(T_1)$, $r_G(u,v)\ge
2+\frac{1\cdot (k-1)}{k}>2$. Let $b$ be the number of neighbors of
$u$ on $C_k$, where $b=0,1$. Now by Eq. (\ref{cycle-vertex}) and inequality (\ref{upp-2}), we have
\begin{eqnarray*} 
Kf_G(u)&=&Kf_{C_k}(u)+\sum_{x \in V (G) \setminus V (C_k)} r_G (u,x) \nonumber \\
       &\ge&\left(k+\frac{k^2-1}{6}\right)+[\left(d_G(u)-b\right)+2(n-k-d_G(u)+b)] \nonumber \\
       &=&b-d_{G}(u)+\frac{1}{6}k^2-k+2n-\frac{1}{6} \nonumber \\
       &\ge&0-\left(n-m + 1 -\left\lfloor \frac{k + 1} {2} \right \rfloor\right)+\frac{1}{6}k^2-k+2n-\frac{1}{6} \nonumber \\
       &>&\frac{1}{6}(k^2-9k+6n+6m-4).
\end{eqnarray*}

Now combining Cases 1 and 2, we have
\[
Kf_G(u) \ge \frac{1}{6}(k^2-9k+6n+6m-4)
\]
with equality if and only if $G \cong U (k,1,n-2m,m-\frac{k+1}{2})$
with odd $k$, and $u$ is the vertex of maximum degree in
$U(k,1,n-2m,m-\frac{k+1}{2})$.  Thus
\begin{eqnarray*}
Kf_G(u) \ge \frac{1}{6}(k^2-9k+6n+6m-4)
&\ge&\frac{1}{6}(5^2-9\cdot 5+6n+6m-4)\\
&=&n+m-4
\end{eqnarray*}
with equalities if and only if $G \cong U(5,1,n-2m,m-3) = U_{n,m}$,
and $u$ is the vertex of maximum degree in $U_{n,m}$.
\end{Proof}

Now we present a stronger version of lemma \ref{degree}.

\begin{Lemma}\label{G-u}
Let $G$ be a unicyclic graph with $n$ vertices and matching number at least $m$,
where $n \ge 6$, $m \ge 3$. For $u\in V(G)$,
\[
Kf_{G}(u)\ge n+m-4
\]
with equality if and only if $G \cong U_{n,m}$, and $u$ is the vertex
of maximum degree in $U_{n,m}$.
\end{Lemma}

\begin{Proof}
Let $G$ be a unicyclic graph with a vertex $u \in V (G)$ such that
\begin{eqnarray} \label{kf-min}
Kf_G (u) = \min \{ Kf_H (x) : H \in \mathbb{U}(n,r), x \in V (H), r \ge m \}.
\end{eqnarray}
Assume that $u \in V(T_1)$. Let $M$ be a maximum matching of $G$.
Suppose that $|T_i| \ge 3$ for some $i$ with $2 \le i \le k$, where
$k$ is the length of the unique cycle of $G$. Then there is
some edge, say $x y$, in $T_i$ outside $M$. Assume that the vertices
$x$ and $u$ lie in the same component of $G - xy$. Let $G_1 = G -xy
+ uy$. Clearly, $M$ is also a matching of $G_1$, and thus $G_1$ has
matching number at least $m$. However, $Kf_{G_1} (u) < Kf_G
(u)$, which is a contradiction. Thus $|T_i| = 1, 2$, i.e., $T_i \cong
P_1$ or $P_2$ for $2 \le i \le k$. By Lemma \ref{degree}, we have
\[
Kf_{G}(u)\ge n+r-4 \ge n+m-4
\]
with equalities if and only if $G \cong U_{n,m}$, and $u$ is the
vertex of maximum degree in $U_{n,m}$.
\end{Proof}

The following result turns out to be of rather important for the
proof of our main results.

\begin{Lemma}  \label{G-x-y}
Let $G\in \mathbb{U}(n,m)$ with a pendent vertex $x$
being adjacent to vertex $y$, and let $z$ be the neighbor of $y$
different from $x$ if $d_G(y)=2$, where $n \ge 6$, $m \ge 3$. Then
\[
Kf(G)-Kf(G-x)\ge 2n+m-6
\]
with equality if and only if $G \cong U_{n,m}$, and $x$ is a pendent neighbor of the vertex of maximum degree in $U_{n,m}$. Moreover, if
$d_G(y)=2$, then
\[
Kf(G)-Kf(G-x-y)\ge5n+2m-19
\]
with equality if and only if $G \cong U_{n,m}$.
\end{Lemma}

\begin{Proof}
Note that $Kf_G(x)-Kf_G(y)=n-2$. Then by Lemma \ref{G-u}, we have
\begin{eqnarray*}
Kf(G)-Kf(G-x)&=&Kf_G(x)\\
&=&Kf_G(y)+n-2\\
       &\ge&(n+m-4)+n-2=2n+m-6
\end{eqnarray*}
with equality if and only if $G \cong U_{n,m}$, and $x$ is a pendent neighbor of the vertex of maximum degree in $U_{n,m}$.

If $d_G(y)=2$, then
$Kf_G(y)-Kf_G(z)=n-4$, and thus by Lemma \ref{G-u}, we have
\begin{eqnarray*}
Kf(G)-Kf(G-x-y)&=&Kf_G(x)+Kf_G(y)-1\\
&=&2Kf_G(z)+3n-11\\
       &\ge&2(n+m-4)+3n-11=5n+2m-19
\end{eqnarray*}
with equality if and only if $G \cong U_{n,m}$.
\end{Proof}

\begin{Lemma}  \label{S-n-k}\textnormal{\cite{YJ}}
Let $G$ be an $n$-vertex unicyclic graph with the unique cycle $C_k$,
where $3\le k\le n-1$. Then
\[
Kf (G) \ge \frac{1}{12}[-k^3+2nk^2-(12n - 13)k+12n^2-14n]
\]
with equality if and only if $G \cong U(k,1,n-k-1,0)$.
\end{Lemma}

\section{Results}


First we consider the minimum Kirchhoff index of unicyclic graphs
with perfect matching.

\begin{Theorem}  \label{min-2m-m}
Among the graphs in $\mathbb{U}(2m,m)$ with $m\ge2$, $C_{2m}$ for
$2\le m\le4$, $U(8,2)$ for $m=5$, $U(8,4)$ for $m=6$, $U(7,7)$ for
$m=7$, and $U_{2m,m}$ for $m\ge8$ are the unique graphs with the
minimum Kirchhoff indices, which are equal to $\frac{1}{6}(4m^3-m)$
for $2\le m\le4$, $81\frac{7}{8}$ for $m=5$, $135\frac{1}{2}$ for
$m=6$, $203$ for $m=7$, and $6m^2-13m+4$ for $m\ge8$.
\end{Theorem}

\begin{Proof}
Recall that $\mathbb{U}(2m,m) = \mathbb{U}_1(m) \cup
\mathbb{U}_2(m) \cup \{ C_{2m} \}$.


\noindent {\bf Case 1.}  $2\le m\le 8$.

%
%
%
%

The case $m=2$ is obvious since $\mathbb{U}(4,2)=\{U (3,1),C_4\}$,
where
$$
Kf(U (3,1)) =6\frac{1}{3}>5=Kf(C_4).
$$

For $3 \le m \le 8$, by Lemmas \ref{U-2-m} and \ref{U-1-m}, the
minimum Kirchhoff index of the graphs in $\mathbb{U}(2m,m)$ is
precisely achieved by some graph of the form $U(k,t,0,j)$, where
$k+t+2j =2m$, $k \ge
3$, $k \ge t \ge 0$ and $j \ge 0$.
In Tables 1--6 corresponding to $m=3,4, \dots, 8$, we
list these graphs and their Kirchhoff indices. We use $(k,t;j)$ to
represent the graph $U(k,t,0,j)$ in these tables. From these tables,
we find that
\begin{enumerate}
\item[(1)]
$U (6,0,0,0) = C_6$ is the unique graph in $\mathbb{U}(6,3)$ with the
minimum Kirchhoff index, which is equal to $17 \frac{1}{2}$$;$

\item[(2)]
$U (8,0,0,0) = C_8$ is the unique graph  in $\mathbb{U}(8,4)$ with the
minimum Kirchhoff index, which is equal to $42$$;$

\item[(3)]
$U (8,2,0,0) = U (8,2)$ is the unique graph in $\mathbb{U}(10,5)$
with the minimum Kirchhoff index, which is equal to $81\frac{7}{8}$$;$

\item[(4)]
$U (8,4,0,0) = U(8,4)$ is the unique graph in $\mathbb{U}(12,6)$
with the minimum Kirchhoff index, which is equal to $135\frac{1}{2}$$;$

\item[(5)]
$U (7,7,0,0) = U (7,7)$ is the unique graph in $\mathbb{U}(14,7)$
with the minimum Kirchhoff index, which is equal to $203$$;$

\item[(6)]
$U (5,1,0,5) = U_{16,8}$ is the unique graph in $\mathbb{U}(16,8)$
with the minimum Kirchhoff index, which is equal to $284$.
\end{enumerate}


\begin{table}[h]
\caption{The Kirchhoff indices of the graphs $U(k,t,0,j)$ in
$\mathbb{U}(6,3)$.}
\bigskip
\renewcommand{\arraystretch}{1.5}
\centering
\begin{tabular}{|c||c|c|c|c|c|c|}
\noalign{\hrule height 1.2pt}
{\text Graphs} & $(3,1;1)$ & $(3,3;0)$ & $(4,0;1)$ & $(4,2;0)$ & $(5,1;0)$ & $(6,0;0)$\\
\hline
{\text Kirchhoff indices} &$24$ & $23$ & $23$ & $20 \frac{3}{4}$ & $19$ & $17 \frac{1}{2}$ \\
\noalign{\hrule height 1.2pt}
\end{tabular}
\end{table}


\begin{table}[!htb]
\small \caption{The Kirchhoff indices of the graphs $U(k,t,0,j)$ in
$\mathbb{U}(8,4)$.}
\bigskip
\renewcommand{\arraystretch}{1.5}
\centering
\begin{tabular}{|c||c|c|c|c|c|c|}
\noalign{\hrule height 1.2pt}
{\text Graphs} & $(3,1;2)$ & $(3,3;1)$ & $(4,0;2)$ & $(4,2;1)$ & $(4,4;0)$ & $(5,1;1)$ \\
\hline
{\text Kirchhoff indices} &$53 \frac{2}{3}$ & $53 \frac{1}{3}$ & $53$ & $50 \frac{1}{4}$ & $48$ & $48$ \\
\noalign{\hrule height 1.2pt}
{\text Graphs} & $(5,3;0)$ & $(6,0;1)$& $(6,2;0)$  & $(7,1;0)$ & $(8,0;0)$ &  \\
\hline
{\text Kirchhoff indices} &$45 \frac{4}{5}$ & $48 \frac{1}{6}$ & $44$  & $43$ & $42$ & \\
\noalign{\hrule height 1.2pt}
\end{tabular}
\end{table}

\begin{table}[!htb]
\small \caption{The Kirchhoff indices of the graphs $U(k,t,0,j)$ in
$\mathbb{U}(10,5)$.}
\bigskip
\renewcommand{\arraystretch}{1.5}
\centering
\begin{tabular}{|c||c|c|c|c|c|c|}
\noalign{\hrule height 1.2pt}
{\text Graphs} & $(3,1;3)$ & $(3,3;2)$ & $(4,0;3)$ & $(4,2;2)$ & $(4,4;1)$ & $(5,1;2)$ \\
\hline
{\text Kirchhoff indices} &$95 \frac{1}{3}$ & $95 \frac{2}{3}$ & $95$ & $91 \frac{3}{4}$ & $91$ & $89$ \\
\noalign{\hrule height 1.2pt}
{\text Graphs} & $(5,3;1)$ & $(5,5;0)$& $(6,0;2)$ & $(6,2;1)$  & $(6,4;0)$ & $(7,1;1)$   \\
\hline
{\text Kirchhoff indices} &$88$ & $85$ & $90 \frac{5}{6}$  & $86 \frac{1}{3}$ & $83 \frac{1}{2}$ & $86$ \\
\noalign{\hrule height 1.2pt}
{\text Graphs} & $(7,3;0)$ & $(8,0;1)$& $(8,2;0)$ & $(9,1;0)$  & $(10,0;0)$ &  \\
\hline
{\text Kirchhoff indices} &$82 \frac{1}{7}$ & $88$ & $81 \frac{7}{8}$  & $82 \frac{1}{3}$ & $82 \frac{1}{2}$ & \\
\noalign{\hrule height 1.2pt}
\end{tabular}
\end{table}

\begin{table}[!htb]
\small \caption{The Kirchhoff indices of the graphs $U(k,t,0,j)$ in
$\mathbb{U}(12,6)$.}
\bigskip
\renewcommand{\arraystretch}{1.5}
\centering
\begin{tabular}{|c||c|c|c|c|c|c|}
\noalign{\hrule height 1.2pt}
{\text Graphs} & $(3,1;4)$ & $(3,3;3)$ & $(4,0;4)$ & $(4,2;3)$ & $(4,4;2)$ & $(5,1;3)$ \\
\hline
{\text Kirchhoff indices} &$149$ & $150$ & $149$ & $145 \frac{1}{4}$ & $146$ & $142$ \\
\noalign{\hrule height 1.2pt}
{\text Graphs} & $(5,3;2)$ & $(5,5;1)$& $(6,0;3)$ & $(6,2;2)$  & $(6,4;1)$ & $(6,6;0)$   \\
\hline
{\text Kirchhoff indices} &$142 \frac{1}{5}$ & $142$ & $145 \frac{1}{2}$  & $140 \frac{2}{3}$ & $140 \frac{1}{6}$ & $136$ \\
\noalign{\hrule height 1.2pt}
{\text Graphs} & $(7,1;2)$ & $(7,3;1)$ & $(7,5;0)$ & $(8,0;2)$ & $(8,2;1)$ & $(8,4;0)$  \\
\hline
{\text Kirchhoff indices} &$141$ & $138 \frac{4}{7}$ & $135 \frac{6} {7}$  & $146$ & $139 \frac{5} {8}$ & $135 \frac{1} {2}$  \\
\noalign{\hrule height 1.2pt}
{\text Graphs} & $(9,1;1)$ & $(9,3;0)$ & $(10,0;1)$ & $(10,2;0)$ & $(11,1;0)$ & $(12,0;0)$  \\
\hline
{\text Kirchhoff indices} &$142$ & $136 \frac{1} {3}$ & $146 \frac{1} {2}$  & $138 \frac{2} {5}$ & $141$ & $143$  \\
\noalign{\hrule height 1.2pt}
\end{tabular}
\end{table}

\begin{table}[!htb]
\small \caption{The Kirchhoff indices of the graphs $U(k,t,0,j)$ in
$\mathbb{U}(14,7)$.}
\bigskip
\renewcommand{\arraystretch}{1.5}
\centering
\begin{tabular}{|c||c|c|c|c|c|c|}
\noalign{\hrule height 1.2pt}
{\text Graphs} & $(3,1;5)$ & $(3,3;4)$ & $(4,0;5)$ & $(4,2;4)$ & $(4,4;3)$ & $(5,1;4)$ \\
\hline
{\text Kirchhoff indices} &$214 \frac{2}{3}$ & $216 \frac{1}{3}$ & $215$ & $210 \frac{3}{4}$ & $213$ & $207$ \\
\noalign{\hrule height 1.2pt}
{\text Graphs} & $(5,3;3)$ & $(5,5;2)$& $(6,0;4)$ & $(6,2;3)$  & $(6,4;2)$ & $(6,6;1)$   \\
\hline
{\text Kirchhoff indices} &$208 \frac{2}{5}$ & $211$ & $212 \frac{1}{6}$  & $207$ & $208 \frac{5}{6}$ & $208 \frac{1}{3}$ \\
\noalign{\hrule height 1.2pt}
{\text Graphs} & $(7,1;3)$ & $(7,3;2)$ & $(7,5;1)$ & $(7,7;0)$ & $(8,0;3)$ & $(8,2;2)$  \\
\hline
{\text Kirchhoff indices} &$208$ & $207$ & $208$  & $203$ & $216$ & $209 \frac{3}{8}$  \\
\noalign{\hrule height 1.2pt}
{\text Graphs} & $(8,4;1)$ & $(8,6;0)$ & $(9,1;2)$ & $(9,3;1)$ & $(9,5;0)$ & $(10,0;2)$  \\
\hline
{\text Kirchhoff indices} &$208$ & $204 \frac{7}{8}$ & $213 \frac{2}{3}$  & $209 \frac{5}{9}$ & $206 \frac{1}{9}$ & $222 \frac{1}{2}$  \\
\noalign{\hrule height 1.2pt}
{\text Graphs} & $(10,2;1)$ & $(10,4;0)$ & $(11,1;1)$ & $(11,3;0)$ & $(12,0;1)$ & $(12,2;0)$  \\
\hline
{\text Kirchhoff indices} &$214 \frac{1}{5}$ & $208 \frac{1}{2}$ & $220$  & $212 \frac{9}{20}$ & $227 \frac{2}{3}$ & $217 \frac{7}{12}$  \\
\noalign{\hrule height 1.2pt}
{\text Graphs} & $(13,1;0)$ & $(14,0;0)$ &  &  & &   \\
\hline
{\text Kirchhoff indices} &$223$ & $227 \frac{1} {2}$ &  & &  &  \\
\noalign{\hrule height 1.2pt}
\end{tabular}
\end{table}

\begin{table}[!htb]
\small \caption{The Kirchhoff indices of the graphs $U(k,t,0,j)$ in
$\mathbb{U}(16,8)$.}
\bigskip
\renewcommand{\arraystretch}{1.5}
\centering
\begin{tabular}{|c||c|c|c|c|c|c|}
\noalign{\hrule height 1.2pt}
{\text Graphs} & $(3,1;6)$ & $(3,3;5)$ & $(4,0;6)$ & $(4,2;5)$ & $(4,4;4)$ & $(5,1;5)$ \\
\hline
{\text Kirchhoff indices} &$292 \frac{1}{3}$ & $293 \frac{2}{3}$ & $293$ & $288 \frac{1}{4}$ & $292$ & $284$ \\
\noalign{\hrule height 1.2pt}
{\text Graphs} & $(5,3;4)$ & $(5,5;3)$& $(6,0;5)$ & $(6,2;4)$  & $(6,4;3)$ & $(6,6;2)$   \\
\hline
{\text Kirchhoff indices} &$286 \frac{3}{5}$ & $292$ & $290 \frac{5}{6}$  & $285 \frac{1}{3}$ & $289 \frac{1}{2}$ & $292 \frac{2}{3}$ \\
\noalign{\hrule height 1.2pt}
{\text Graphs} & $(7,1;4)$ & $(7,3;3)$ & $(7,5;2)$ & $(7,7;1)$ & $(8,0;4)$ & $(8,2;3)$  \\
\hline
{\text Kirchhoff indices} &$287$ & $287 \frac{3}{7}$ & $292 \frac{1}{7}$  & $292$ & $298$ & $291 \frac{1}{8}$  \\
\noalign{\hrule height 1.2pt}
{\text Graphs} & $(8,4;2)$ & $(8,6;1)$ & $(8,8;0)$ & $(9,1;3)$ & $(9,3;2)$ & $(9,5;1)$ \\
\hline
{\text Kirchhoff indices} & $292 \frac{1}{2}$ & $294 \frac{1}{8}$ & $288$  & $297 \frac{1}{3}$ & $294 \frac{7}{9}$ & $295 \frac{5}{9}$  \\
\noalign{\hrule height 1.2pt}
{\text Graphs} & $(9,7;0)$ & $(10,0;3)$ & $(10,2;2)$ & $(10,4;1)$ & $(10,6;0)$ & $(11,1;2)$ \\
\hline
{\text Kirchhoff indices} &$292 \frac {5} {9}$ & $310 \frac{1}{2}$ & $302$  & $299 \frac{3}{10}$ & $296$ & $311$  \\
\noalign{\hrule height 1.2pt}
{\text Graphs} & $(11,3;1)$ & $(11,5;0)$ & $(12,0;2)$ & $(12,2;1)$ & $(12,4;0)$ & $(13,1;1)$    \\
\hline
{\text Kirchhoff indices} &$305 \frac{1}{11}$ & $300 \frac {5} {11}$ & $324 \frac{1}{3}$ & $314 \frac {1} {12}$ & $306 \frac{2}{3}$ & $324$ \\
\noalign{\hrule height 1.2pt}
{\text Graphs} & $(13,3;0)$ & $(14,0;1)$ & $(14,2;0)$ & $(15,1;0)$ & $(16,0;0)$ &    \\
\hline
{\text Kirchhoff indices} &$314 \frac{7}{13}$ & $335 \frac {1} {2}$ & $323 \frac{3}{7}$ & $332 \frac{1}{3}$ & $340$ & \\
\noalign{\hrule height 1.2pt}
\end{tabular}
\end{table}

%
%

\noindent {\bf Case 2.} $m\ge 9$.

We prove the result by induction on $m$. Suppose that the result
holds for all the graphs in $\mathbb{U}(2m-2,m-1)$.
Let $G\in \mathbb{U}(2m,m)$.


If $G \cong C_{2m}$, then by Eqs. (\ref{cycle-all}) and (\ref{eq-Unm}), we have
$$
Kf(G)=\frac{1}{6}(4m^3-m)>6m^2-13m+4 = Kf (U_{2m,m}).
$$

Suppose that $G\in\mathbb{U}_1(m)$. Recall that $G$ is a graph
of maximum degree three obtainable by attaching some pendent
vertices to a cycle $C_k$, where $m\le k\le2m-1$. If $k = m, m + 1,
m + 2$, then there are, respectively, $m, m - 1, m - 2$ pendent
vertices in $G$ outside the cycle $C_k$, and thus by Lemma \ref{U-2-m-a} (i), we have
\begin{eqnarray*}
Kf (G)&\ge&\begin{cases}
\frac{1}{3}(m^3+6m^2-4m) & \text{if $k = m$}\\
\frac{1}{3}\left(m^3+7m^2-11m+6-\frac{3}{m+1}\right) & \text{if $k = m + 1$}\\
\frac{1}{3}\left(m^3+8m^2-20m+30-\frac{60}{m+2}\right) & \text{if $k = m + 2$}
            \end{cases}\\
&>& 6m^2-13m+4 = Kf (U_{2m,m}).
\end{eqnarray*}
If $m+3\le k\le 2m-1$, then by Lemma \ref{S-n-k}, $Kf(G) \ge
\frac{1}{12} h(k)$, where $h(k)= -k^3+4mk^2 - (24m-13)k+48m^2-28m$.
Clearly, $h'(k)=-3k^2+8mk-24m+13$. Note that $h'(m+3)=5m^2-18m-14>0$ and
$h'(2m-1)=4m^2-20m+10>0$. This implies that  $h'(k)>0$ for $m+3\le
k\le 2m-1$, i.e., $h(k)$ is increasing for $k$ with $m+3\le k\le
2m-1$. Thus
\[
Kf(G)\ge \frac{1}{12}h(m+3)=\frac{1}{4}(m^3+13m^2-26m+4)>6m^2-13m+4 = Kf (U_{2m,m}).
\]

Now suppose that $G\in \mathbb{U}_2(m)$.
Denote by $x$ a pendent vertex
in $G$ whose
unique neighbor $y$ is of degree two,
and $z$ the neighbor of
$y$ different from $x$ in $G$. Obviously, $xy\in M$. Then
$G-x-y\in\mathbb{U}(2m-2,m-1)$,
and thus by Lemma \ref{G-x-y} and
the induction hypothesis, we have
\begin{eqnarray*}
Kf(G)&\ge&Kf(G-x-y)+12m-19\\
     &\ge& Kf (U_{2m - 2,m - 1})+12m-19\\
     &=&6m^2-13m+4= Kf (U_{2m,m})
\end{eqnarray*}
with equalities if and only if $G \cong U_{2m,m}$.


Then the result for $m \ge 9$ follows easily.
\end{Proof}

The remainder of the paper will focus on the minimum Kirchhoff index among the graphs in $\mathbb{U}(n,m)$, where $n > 2m$ and $m \ge 3$.

\begin{Lemma}\textnormal{\cite{CT}}\label{l4}
Let $G\in\mathbb {U}(n,m) \setminus \{C_n\}$, where $n> 2m$, $m \ge 3$. Then
there is a maximum  matching $M$ and a pendent vertex $u$ of $G$
such that $u$ is not  $M$-saturated.
\end{Lemma}


For a given graph $G\in\mathbb {U}(n,m) \setminus \{C_n\}$, where $n> 2m$, $m \ge 3$,
by Lemma \ref{l4}, there
is a maximum matching $M$ and a pendent vertex which is not
$M$-saturated, after deleting this pendent vertex, we can get a graph in $\mathbb {U}(n-1,m)$. Repeating
the process until it is exhausted, the resulting graph is denoted by $G_0$.
Note that $G_0 \in \mathbb{U}(2m,m)$. Let $n_0 = |G_0|$.
Furthermore, for the vertex $u \in V (G_0)$ satisfying $Kf_{G_0} (u)$ is minimum, let $G_0^*$ be the graph obtained from $G_0$ by attaching $n - n_0$ pendent vertices to $u$.

\begin{Lemma}\label{adddddd}
Let $G \in \mathbb{U}(n,m) \setminus \{C_n\}$, where $n> 2m$, $m \ge 3$.
\begin{enumerate}
\item[$(i)$]
Then $Kf (G) \ge Kf (G_0^*)$. In particular, if $u$ is the unique vertex in $G_0$ such that $Kf_{G_0} (u)$ is minimum, then $Kf (G) \ge Kf (G_0^*)$ with equality if and only if $G \cong G_0^*$.
\item[$(ii)$]
If $G_0 \not\cong U_{n_0,m}$ and $Kf (G_0) \ge Kf (U_{n_0,m})$, then $Kf (G) > Kf (U_{n,m})$.
\end{enumerate}
\end{Lemma}

\begin{Proof}
%
Similar to the proof of Lemma \ref{U-1-m}, and by Lemma \ref{formula} repeatedly,
(i) follows easily. On the other hand, by Lemmas \ref{formula} and \ref{G-u},
$Kf (G_0^*) > Kf (U_{n,m})$ follows from the hypothesis that $G_0 \not\cong U_{n_0,m}$ and $Kf (G_0) \ge Kf (U_{n_0,m})$.
Now together with $Kf (G) \ge Kf (G_0^*)$, we can get $Kf (G) > Kf (U_{n,m})$.
\end{Proof}

%
%
%
%
%
%
%
%
%
%
%
%

The following lemma reveal the possible graph with the minimum Kirchhoff index among the graphs in $\mathbb{U}(n,m) \setminus \{C_n\}$, where $n > 2m$ and $3\le m\le7$.

\begin{Lemma} \label{aadd}
If $G$ is a graph in $\mathbb{U}(n,m) \setminus \{C_n\}$ with the minimum Kirchhoff index,
where $n > 2m$ and $3\le m\le7$, then $G_0 \cong U (k,t,0,j)$ with $k+t +2j=n_0$, $k \ge
3$, $k \ge t \ge 0$ and $j \ge 0$.
\end{Lemma}

\begin{Proof}
Let $G \in \mathbb{U}(n,m) \setminus \{C_n\}$. Denote by $k$ the length of the unique cycle of $G$. Suppose that there are $t$ pendent vertices of $G$ whose unique neighbors are all on the unique cycle of $G$. Note that $G_0 \in \mathbb{U}(2m,m)$, i.e., $G_0 \in \mathbb{U}_1 (m) \cup \mathbb{U}_2 (m)$.

\noindent {\bf Case 1.} $G_0 \in \mathbb{U}_1 (m)$.

First, by Lemma \ref{adddddd} (i), we have $Kf (G) \ge Kf (G_0^*)$. Next, by Lemma \ref{U-2-m}, we have $Kf (G_0) \ge Kf (U(k,t))$ with equality if and only if $G_0 \cong U(k,t))$, and thus by Lemma \ref{U-2-m-a} (ii) and Lemma \ref{formula}, we have $Kf (G_0^*) \ge Kf (U(k,t,i,0))$ with equality if and only if $G_0^* \cong U(k,t,i,0))$, where $k + t + i = n$ and $i \ge 1$.
Now it follows that $Kf (G) \ge Kf (U(k,t,i,0))$ with equality if and only if $G \cong U(k,t,i,0)$.

\noindent {\bf Case 2.} $G_0 \in \mathbb{U}_2 (m)$.

Recall that,
starting from $G_0$, deleting the pendent
$P_2$'s repeatedly, until there is no pendent
$P_2$,
the resulting graph is denoted by $\bar{G_0}$. Let $\bar{n_0} = |\bar{G_0}|$. Clearly, $\bar{G_0} \in \mathbb{U}_1 (\frac{\bar{n_0}} {2}) \cup \{C_{\bar{n_0}}\}$.

Suppose that $u$ is a vertex in $\bar{G_0}$ satisfying $Kf_{\bar{G_0}} (u)$ is minimum, and let $H$ be the graph obtained from $\bar{G_0}$ by attaching $i$ pendent vertices and $j$ paths on two vertices to $u$.

Similar to the proof of Lemma \ref{U-1-m}, and by Lemma \ref{formula} repeatedly,
$Kf (G) \ge Kf (H)$ follows easily.

On the other hand, recall that $\bar{G_0} \in \mathbb{U}_1 (\frac{\bar{n_0}} {2}) \cup \{C_{\bar{n_0}}\}$, by Lemma \ref{U-2-m}, we have $Kf (\bar{G_0}) \ge Kf (U(k,t))$ with equality if and only if $\bar{G_0} \cong U(k,t))$, where $k + t = \bar{n_0}$, and thus by Lemma \ref{U-2-m-a} (ii) and Lemma \ref{formula}, we have $Kf (H) \ge Kf (U (k,t,i,j))$ with equality if and only if $H \cong U (k,t,i,j)$, where $k+t+i+2j =n$.

Now it follows that $Kf (G) \ge Kf (U (k,t,i,j))$ with equality if and only if $G \cong U (k,t,i,j)$.

Combining Cases 1 and 2, we have $Kf (G) \ge Kf (U (k,t,i,j))$ with equality if and only if $G \cong U (k,t,i,j)$, and $G \cong U (k,t,i,j)$ implies that $G_0 \cong U (k,t,0,j)$.
\end{Proof}

Now we determine the minimum Kirchhoff index among the unicyclic graphs
with given matching number.

\begin{Theorem}  \label{min-n-m}
Among the graphs in $\mathbb{U}(n,m)$ with $2\le
m\le\lfloor\frac{n}{2}\rfloor$,
\begin{enumerate}
\item[$(i)$]
for $m=2$, $C_n$ for $n=4,5$, $U(4,1,n-5,0)$ for $6\le
n\le11$, $U(3,1,8,0)$ and $U(4,1,7,0)$ for $n=12$, and
$U(3,1,n-4,0)$ for $n\ge13$ are the unique graphs with the minimum
Kirchhoff indices, which are equal to $\frac{n^3-n}{12}$ for $n=4,5$,
$\frac{1}{2}(2n^2-5n-2)$ for $6\le n\le11$, $113$ for $n=12$, and
$\frac{1}{3}(3n^2-8n+3)$ for $n\ge13$$;$
\item[$(ii)$]
for $m=3$, $C_n$ for $n=6,7$  and $U_{n,3}$ for $n\ge8$ are
the unique graphs with the minimum Kirchhoff indices, which are equal to
$\frac{n^3-n}{12}$ for $n=6,7$ and
$n^2-2n-5$ for $n\ge8$$;$
\item[$(iii)$]
for $m=4$, $C_8$ for $n=8$, $U(7,1,1,0)$ and $C_9$ for
$n=9$, $U(7,1,2,0)$ for $n=10$, $U(6,2,3,0)$ and $U(7,1,3,0)$ for
$n=11$, $U(6,2,n-8,0)$ for $n=12,13$, $U_{14,4}$ and $U(6,2,6,0)$
for $n=14$, and $U_{n,m}$ for $n\ge15$ are the unique graphs with the
minimum Kirchhoff indices, which are equal to $42$ for $n=8$, $60$
for $n=9$, $79$ for $n=10$, $100$ for $n=11$,
$\frac{1}{3}(3n^2-n-52)$ for $n=12,13$, $174$ for $n=14$, and
$n^2-n-8$ for $n\ge15$$;$
\item[$(iv)$]
for $m=5$, $U(8,2)$ for $n=10$, $U(7,3,n-10,0)$ for $11\le
n\le13$, $U_{14,5}$, $U(6,2,4,1)$ and $U(7,3,4,0)$ for $n=14$, and
$U_{n,5}$ for $n\ge15$ are the unique graphs with the minimum Kirchhoff
indices, which are equal to $81\frac{7}{8}$ for $n=10$,
$\frac{1}{7}(7n^2+12n-245)$ for $11\le n\le13$, $185$ for $n=14$,
and $n^2-11$ for $n\ge15$$;$
\item[$(v)$]
for $m=6$, $U(8,4,n-12,0)$ for $n=12,13$, $U_{14,6}$,
$U(6,2,2,2)$  and $U(7,3,2,1)$ for $n=14$, and $U_{n,6}$ for
$n\ge15$ are the unique graphs with the minimum Kirchhoff indices, which
are equal to $\frac{1}{4}(4n^2+19n-262)$ for $n=12,13$, $196$ for
$n=14$, and
$n^2+n-14$ for $n\ge15$$;$
\item[$(vi)$]
for $m=7$, $U(7,7)$ for $n=14$  and $U_{n,7}$ for $n\ge15$
are the unique graphs with the minimum Kirchhoff indices, which are equal
to $203$ for $n=14$ and
$n^2+2n-17$ for $n\ge15$$;$
\item[$(vii)$]
for $m\ge8$, $U_{n,m}$ for $n\ge16$ is the unique graph with the
minimum Kirchhoff index, which is equal to $n^2+nm-5n-3m+4$.
\end{enumerate}
\end{Theorem}

\begin{Proof}
The result for $n=2m$ follows from Theorem \ref{min-2m-m}. Suppose
that $n>2m$. Let $G \in \mathbb{U}(n,m)$.





\noindent {\bf Case 1.} $m=2$.

Clearly, the girth of $G$ is $3,4,5$ for
$n=5$, and $3,4$ for $n\ge 6$.  Then by Lemma \ref{S-n-k}, we have
\begin{eqnarray*}
Kf(G)&\ge&\min\{Kf(U(3,1,1,0)),Kf(U(4,1)),Kf(C_5)\}\\
     &=&\min\left\{12\frac{2}{3},11\frac{1}{2}, 10\right\}=10
\end{eqnarray*}
for $n=5$, and
\begin{eqnarray*}
Kf(G)&\ge&\min\{Kf(U(3,1,n-4,0)),Kf(U(4,1,n-5,0))\}\\
     &=&\min\left\{\frac{1}{3}(3n^2-8n+3),\frac{1}{2}(2n^2-5n-2)\right\}
\end{eqnarray*}
for $n\ge 6$. Thus $C_5$ for $n=5$, $U(4,1,n-5,0)$ for $6\le n\le11$,
$U(3,1,8,0)$ and $U(4,1,7,0)$ for $n=12$, and $U(3,1,n-4,0)$ for
$n\ge13$ are the unique graphs in $\mathbb{U}(n,2)$ with the minimum
Kirchhoff indices.



\noindent {\bf Case 2.} $m=3$.

If $G \cong C_n$, then $n = 7$, and by Eq. (\ref{cycle-all}), we have $Kf (G) = 28$.

Suppose that $G \not\cong C_n$.
If $G_0 \cong U_{6,3}$, then by Lemma \ref{adddddd} (i), we have $Kf (G) \ge Kf (U_{n,3})$ with equality if and only if $G \cong U_{n,3}$.
Suppose that $G_0 \not\cong U_{6,3}$.
If $Kf (G_0) \ge Kf (U_{6,3})$, then by Lemma \ref{adddddd} (ii), we have $Kf (G) > Kf (U_{n,3})$.
If $Kf (G_0) < Kf (U_{6,3})$, then by Lemma \ref{aadd} and Table 1, we assume that $G_0 = C_6$, and thus by Lemma \ref{adddddd} (i), we have $Kf (G) \ge Kf (U(6,1,n-7,0))$ with equality if and only if $G \cong U(6,1,n-7,0)$.

Therefore for $n = 7$,
$$
Kf(G) \ge \min\{Kf (U_{7,3}), Kf (U(6,1)),
Kf(C_7)\}=\min\left\{30,29\frac{1}{3},28\right\}=28
$$
with equality if and only if $G \cong C_7$,
and for $n \ge 8$,
\begin{eqnarray*}
Kf(G) &\ge& \min\{Kf (U_{n,3}), Kf (U(6,1,n-7,0))\} \\
&=&\min\left\{n^2-2n-5,n^2-\frac{7}{6}n-\frac{23}{2}\right\}=n^2-2n-5
\end{eqnarray*}
with equality if and only if $G \cong U_{n,3}$.

\noindent {\bf Case 3.} $m=4$.

If $G \cong C_n$, then $n = 9$, and by Eq. (\ref{cycle-all}), we have $Kf (G) = 60$.

Suppose that $G \not\cong C_n$.
If $G_0 \cong U_{8,4}$, then by Lemma \ref{adddddd} (i), we have $Kf (G) \ge Kf (U_{n,4})$ with equality if and only if $G \cong U_{n,4}$.
Suppose that $G_0 \not\cong U_{8,4}$.
If $Kf (G_0) \ge Kf (U_{8,4})$, then by Lemma \ref{adddddd} (ii), we have $Kf (G) > Kf (U_{n,4})$.
If $Kf (G_0) < Kf (U_{8,4})$, then by Lemma \ref{aadd} and Table 2, we assume that $G_0 = U(5,3,0,0),
U(6,2,0,0), U(7,1,0,0)$ or $U(8,0,0,0)$, and thus by Lemma \ref{adddddd} (i), we have
\begin{eqnarray*}
Kf (G) &\ge& \min \{ Kf (U(5,3,n-8,0)), Kf (U(6,2,n-8,0)), \\
&&Kf (U(7,1,n-8,0)), Kf (U(8,0,n-8,0)) \}.
\end{eqnarray*}

Now the result for $m = 4$ follows from Table 7 easily.

\begin{table}[h]
\small \caption{The graphs in $\mathbb{U}(n,4)$ and their Kirchhoff indices.}
\bigskip
\renewcommand{\arraystretch}{1.5}
\centering
\begin{tabular}{|c|cccccccc|}
\noalign{\hrule height 1.2pt}
\raisebox{-1.4ex}[0pt]
{Graphs} & \multicolumn{8}{c|}{Kirchhoff indices}\\
\cline{2-9} & \makebox[1.8cm][c]{$n$} & \makebox[0.8cm][c]{$9$} &
\makebox[0.9cm][c]{$10$} & \makebox[0.9cm][c]{$11$} &
\makebox[0.9cm][c]{$12$} & \makebox[0.8cm][c]{$13$} &
\makebox[0.9cm][c]{$14$} &
\makebox[0.9cm][c]{$15$}\\
\noalign{\hrule height 1.2pt}
$U_{n,4}$ & $n^2-n-8$ & $64$ & $82$ & $102$ & $124$ & $148$ &  ${\bf 174}$ & ${\bf 202}$\\
\hline $U(5,3,n-8,0)$ & $n^2-\frac{2}{5}n-15$ & $62\frac{2}{5}$ & $81$ & $101\frac{3}{5}$ & $124\frac{1}{5}$ & $148\frac{4}{5}$ &  $175\frac{2}{5}$ & $204$\\
\hline $U(6,2,n-8,0)$ & $n^2-\frac{1}{3}n-\frac{52}{3}$ & $60\frac{2}{3}$ & $79\frac{1}{3}$ & ${\bf 100}$ & ${\bf 122\frac{2}{3}}$ & ${\bf 147\frac{1}{3}}$ &  ${\bf 174}$ & $202\frac{2}{3}$\\
\hline $U(7,1,n-8,0)$ & $n^2-21$ & ${\bf 60}$& ${\bf 79}$ & ${\bf 100}$ & $123$ & $148$ & $175$ &  $204$ \\
\hline $U(8,0,n-8,0)$ & $n^2+\frac{3}{2}n-34$ &  $60\frac{1}{2}$& $81$ & $103\frac{1}{2}$ & $128$ & $154\frac{1}{2}$ & $183$ &  $213\frac{1}{2}$ \\
\hline $C_9$ &  &  ${\bf 60}$&  &  &  &  &  &   \\
\noalign{\hrule height 1.2pt}
\end{tabular}
\end{table}

\noindent {\bf Case 4.} $m=5$.

If $G \cong C_n$, then $n = 11$, and by Eq. (\ref{cycle-all}), we have $Kf (G) =110$.

Suppose that $G \not\cong C_n$.
If $G_0 \cong U_{10,5}$, then by Lemma \ref{adddddd} (i), we have $Kf (G) \ge Kf (U_{n,5})$ with equality if and only if $G \cong U_{n,5}$.
Suppose that $G_0 \not\cong U_{10,5}$.
If $Kf (G_0) \ge Kf (U_{10,5})$, then by Lemma \ref{adddddd} (ii), we have $Kf (G) > Kf (U_{n,5})$.
If $Kf (G_0) < Kf (U_{10,5})$, then by Lemma \ref{aadd} and Table 3, we assume that $G_0 =U(5,3,0,1)$, $U(5,5,0,0)$,
$U(6,2,0,1)$, $U(6,4,0,0)$, $U(7,1,0,1)$, $U(7,3,0,0)$, $U(8,0,0,1)$, $U(8,2,0,0)$, $U(9,1,0,0)$ or $U(10,0,0,0)$, and thus by Lemma \ref{adddddd} (i), we have
\begin{eqnarray*}
Kf (G) &\ge& \min \{ Kf(U(5,3,n-10,1)), Kf (U(5,5,n-10,0)), Kf(U(6,2,n-10,1)),\\
       &&Kf (U(6,4,n-10,0)), Kf(U(7,1,n-10,1)), Kf (U(7,3,n-10,0)), \\
       &&Kf(U(8,0,n-10,1)), Kf (U(8,2,n-10,0)), Kf (U(9,1,n-10,0)), \\
       &&Kf (U(10,0,n-10,0))\}.
\end{eqnarray*}

Now the result for $m = 5$ follows from Table 8 easily.

\begin{table}[h]
\small \caption{The graphs in $\mathbb{U}(n,5)$ and their Kirchhoff indices.}
\bigskip
\renewcommand{\arraystretch}{1.5}
\centering
\begin{tabular}{|c|cccccc|}
\noalign{\hrule height 1.2pt}
\raisebox{-1.4ex}[0pt]
{Graphs} & \multicolumn{6}{c|}{Kirchhoff indices}\\
\cline{2-7} &
\makebox[1.8cm][c]{$n$} & \makebox[0.8cm][c]{$11$} &
\makebox[0.9cm][c]{$12$} & \makebox[0.9cm][c]{$13$} &
\makebox[0.9cm][c]{$14$} & \makebox[0.8cm][c]{$15$} \\
\noalign{\hrule height 1.2pt}
$U_{n,5}$ & $n^2-11$ & $110$ & $133$ &  $158$ &  ${\bf 185}$ & ${\bf 214}$\\
\hline $U(5,3,n-10,1)$ & $n^2+\frac{3}{5}n-18$  & $109\frac{3}{5}$ & $133\frac{1}{5}$ & $158\frac{4}{5}$ & $186\frac{2}{5}$ & $216$\\
\hline $U(5,5,n-10,0)$ & $n^2+2n-35$  & $108$ & $133$ & $160$ & $189$ & $220$\\
\hline $U(6,2,n-10,1)$ & $n^2+\frac{2}{3}n-\frac{61}{3}$  & $108$ & $131\frac{2}{3}$ & $157\frac{1}{3}$ & $\bf{185}$ & $214\frac{2}{3}$\\
\hline $U(6,4,n-10,0)$ & $n^2+\frac{11}{6}n-\frac{209}{6}$  & $106\frac{1}{3}$ & $131\frac{1}{6}$ & $158$ & $186\frac{5}{6}$ & $217\frac{2}{3}$\\
\hline $U(7,1,n-10,1)$ & $n^2+n-24$  & $108$ & $132$ & $158$ & $186$ & $216$\\
\hline $U(7,3,n-10,0)$ & $n^2+\frac{12}{7}n-\frac{245}{7}$  & $\bf{104\frac{6}{7}}$ & $\bf{129\frac{4}{7}}$ & $\bf{156\frac{2}{7}}$ & $\bf{185}$ & $215\frac{5}{7}$\\
\hline $U(8,0,n-10,1)$ & $n^2+\frac{5}{2}n-37$  & $111\frac{1}{2}$ & $137$ & $164\frac{1}{2}$ & $194$ & $225\frac{1}{2}$\\
\hline $U(8,2,n-10,0)$ & $n^2+\frac{19}{8}n-\frac{335}{8}$  & $105\frac{1}{4}$ & $130\frac{5}{8}$ & $158$ & $187\frac{3}{8}$ & $218\frac{3}{4}$\\
\hline $U(9,1,n-10,0)$ & $n^2+\frac{10}{3}n-51$  & $106\frac{2}{3}$ & $133$ & $161\frac{1}{3}$ & $191\frac{2}{3}$ & $224$\\
\hline $U(10,0,n-10,0)$ & $n^2+\frac{11}{2}n-\frac{145}{2}$  & $109$ & $137\frac{1}{2}$ & $168$ & $200\frac{1}{2}$ & $235$\\
\hline $C_{11}$ &  & $110$&  &  &  &    \\
\noalign{\hrule height 1.2pt}
\end{tabular}
\end{table}

\noindent {\bf Case 5.} $m=6$.

If $G \cong C_n$, then $n = 13$, and by Eq. (\ref{cycle-all}), we have $Kf (G) =182$.

Suppose that $G \not\cong C_n$.
If $G_0 \cong U_{12,6}$, then by Lemma \ref{adddddd} (i), we have $Kf (G) \ge Kf (U_{n,6})$ with equality if and only if $G \cong U_{n,6}$.
Suppose that $G_0 \not\cong U_{12,6}$.
If $Kf (G_0) \ge Kf (U_{12,6})$, then by Lemma \ref{adddddd} (ii), we have $Kf (G) > Kf (U_{n,6})$.
If $Kf (G_0) < Kf (U_{12,6})$, then by Lemma \ref{aadd} and Table 4, we assume that $G_0 =U(6,2,0,2)$, $U(6,4,0,1)$, $U(6,6,0,0)$, $U(7,1,0,2)$, $U(7,3,0,1)$, $U(7,5,0,0)$, $U(8,2,0,1)$, $U(8,4,0,0)$, $U(9,3,0,0)$, $U(10,2,0,0)$ or $U(11,1,0,0)$, and thus by Lemma \ref{adddddd} (i), we have
\begin{eqnarray*}
Kf (G) &\ge& \min \{ Kf(U(6,2,n-12,2)), Kf (U(6,4,n-12,1)), Kf(U(6,6,n-12,0)),\\
       && Kf(U(7,1,n-12,2)), Kf (U(7,3,n-12,1)), Kf(U(7,5,n-12,0)),\\
       &&Kf(U(8,2,n-12,1)), Kf (U(8,4,n-12,0)), Kf (U(9,3,n-12,0)), \\
       &&Kf (U(10,2,n-12,0)), Kf(U(11,1,n-12,0))\}.
\end{eqnarray*}

Now the result for $m = 6$ follows from Table 9 easily.

\begin{table}[h]
\small \caption{The graphs in $\mathbb{U}(n,6)$ and their Kirchhoff indices.}
\bigskip
\renewcommand{\arraystretch}{1.5}
\centering
\begin{tabular}{|c|cccc|}
\noalign{\hrule height 1.2pt}
\raisebox{-1.4ex}[0pt]
{Graphs} & \multicolumn{4}{c|}{Kirchhoff indices}\\
\cline{2-5} &
\makebox[1.8cm][c]{$n$} & \makebox[0.8cm][c]{$13$} &
\makebox[0.9cm][c]{$14$} & \makebox[0.9cm][c]{$15$}  \\
\noalign{\hrule height 1.2pt}
$U_{n,6}$ & $n^2+n-14$ &  $168$ &  ${\bf 196}$ & ${\bf 226}$\\
\hline $U(6,2,n-12,2)$ & $n^2+\frac{5}{3}n-\frac{70}{3}$  & $167\frac{1}{3}$ & $\bf{196}$ & $226\frac{2}{3}$ \\
\hline $U(6,4,n-12,1)$ & $n^2+\frac{17}{6}n-\frac{227}{6}$  & $168$ & $197\frac{5}{6}$ & $229\frac{2}{3}$ \\
\hline $U(6,6,n-12,0)$ & $n^2+\frac{14}{3}n-64$  & $165\frac{2}{3}$ & $197\frac{1}{3}$ & $231$ \\
\hline $U(7,1,n-12,2)$ & $n^2+2n-27$  & $168$ & $197$ & $228$ \\
\hline $U(7,3,n-12,1)$ & $n^2+\frac{19}{7}n-38$  & $166\frac{2}{7}$ & $\bf{196}$ & $227\frac{5}{7}$ \\
\hline $U(7,5,n-12,0)$ & $n^2+\frac{32}{7}n-63$  & $165\frac{3}{7}$  &$197$ & $230\frac{4}{7}$ \\
\hline $U(8,2,n-12,1)$ & $n^2+\frac{27}{8}n-\frac{359}{8}$  & $168$  &$198\frac{3}{8}$ & $230\frac{3}{4}$ \\
\hline $U(8,4,n-12,0)$ & $n^2+\frac{19}{4}n-\frac{131}{2}$  & $\bf{165\frac{1}{4}}$  &$197$ & $230\frac{3}{4}$ \\
\hline $U(9,3,n-12,0)$ & $n^2+\frac{46}{9}n-69$  & $166\frac{4}{9}$  &$198\frac{5}{9}$ & $232\frac{2}{3}$ \\
\hline $U(10,2,n-12,0)$ & $n^2+\frac{32}{5}n-\frac{412}{5}$  & $169\frac{4}{5}$  &$203\frac{1}{5}$ & $238\frac{3}{5}$ \\
\hline $U(11,1,n-12,0)$ & $n^2+8n-99$  & $174$  &$209$ & $246$ \\
\hline $C_{13}$ &  & $182$ & & \\
\noalign{\hrule height 1.2pt}
\end{tabular}
\end{table}

\noindent {\bf Case 6.} $m=7$.

If $G \cong C_n$, then $n = 15$, and by Eq. (\ref{cycle-all}), we have $Kf (G)=280$.

Suppose that $G \not\cong C_n$.
If $G_0 \cong U_{14,7}$, then by Lemma \ref{adddddd} (i), we have $Kf (G) \ge Kf (U_{n,7})$ with equality if and only if $G \cong U_{n,7}$.
Suppose that $G_0 \not\cong U_{14,7}$.
If $Kf (G_0) \ge Kf (U_{14,7})$, then by Lemma \ref{adddddd} (ii), we have $Kf (G) > Kf (U_{n,7})$.
If $Kf (G_0) < Kf (U_{14,7})$, then by Lemma \ref{aadd} and Table 5, we assume that $G_0 =U(7,7,0,0)$,  $U(8,6,0,0)$ or $U(9,5,0,0)$, and thus by Lemma \ref{adddddd} (i), we have
\begin{eqnarray*}
Kf (G) &\ge& \min \{ Kf(U(7,7,n-14,0)), Kf (U(8,6,n-14,0)), Kf(U(9,5,n-14,0))\}.
\end{eqnarray*}

Now the result for $m = 7$ follows from Table 10 easily.

\begin{table}[h]
\small \caption{The graphs in $\mathbb{U}(n,7)$ and their Kirchhoff indices.}
\bigskip
\renewcommand{\arraystretch}{1.5}
\centering
\begin{tabular}{|c|cc|}
\noalign{\hrule height 1.2pt}
\raisebox{-1.4ex}[0pt]
{Graphs} & \multicolumn{2}{c|}{Kirchhoff indices}\\
\cline{2-3} & \makebox[1.8cm][c]{$n$} &
\makebox[0.9cm][c]{$15$}
\\
\noalign{\hrule height 1.2pt}
$U_{n,7}$ & $n^2+2n-17$ &  ${\bf 238}$\\
\hline $U(7,7,n-14,0)$ & $n^2+8n-105$ &  $240$\\
\hline $U(8,6,n-14,0)$ & $n^2+\frac{65}{8}n-\frac{839}{8}$ &  $242$\\
\hline $U(9,5,n-14,0)$ & $n^2+\frac{74}{9}n-105$ &  $243\frac{1}{3}$\\
\hline $C_{15}$ &  &  $280$   \\
\noalign{\hrule height 1.2pt}
\end{tabular}
\end{table}

\noindent {\bf Case 7.} $m\ge 8$.

If $G \cong C_n$, then $n = 2m + 1$, and
by Eqs. (\ref{cycle-all}) and (\ref{eq-Unm}), we have
\[
Kf(C_{2m+1}) = \frac{1}{3}(2m^3+3m^2+m) > 6m^2 - 8m = Kf (U_{2m + 1,
m}).
\]
Suppose that $G \not\cong C_n$. By Theorem \ref{min-2m-m}, $Kf (G_0) \ge Kf (U_{n_0,m})$. Furthermore, if $G_0 \cong U_{n_0,m}$, then by Lemma \ref{adddddd} (i), we have $Kf (G) \ge Kf (U_{n,m})$ with equality if and only if $G \cong U_{n,m}$, and if $G_0 \not\cong U_{n_0,m}$, then by Lemma \ref{adddddd} (ii), we have $Kf (G) > Kf (U_{n,m})$. Then the result for $m \ge 8$ follows easily.
\end{Proof}

\vspace{3mm}

\noindent {\bf Acknowledgement.}  This work was supported by the Research Fund for the Doctoral Program of Higher Education of China (No.~20124407110002), the Science Foundation of Hebei Normal University (no.~L2012Q02) and the Guangdong Provincial National Science Foundation of China (no.~S2011010005539).


\begin{thebibliography}{99}


\bibitem{BBLK}
D. Bonchev, A.T. Balaban, X. Liu, D.J. Klein, Molecular cyclicity
and centricity of polycyclic graphs, {\it Int. J. Quantum Chem.\/} {\bf 50} (1994)
1--20.

\bibitem{CT} A. Chang, F. Tian,
On the spectral radius of unicyclic graphs with perfect matching,
{\it Linear Algebra Appl.\/} {\bf 370} (2003) 237--250.

\bibitem{DZ} Z. Du, B. Zhou, Minimum Wiener indices of trees and
unicyclic graphs of given matching number, {\it MATCH Commun. Math.
Comput. Chem.\/} {\bf 63} (2010) 101-112.


\bibitem{DEG}
A.A. Dobrynin, R. Entringer, I. Gutman, Wiener index of trees:
Theory and applications, {\it Acta Appl. Math.\/} {\bf 66} (2001)
211--249.


\bibitem{DGKZ} A.A. Dobrynin, I. Gutman, S. Klav\v{z}ar, P. \v{Z}igert,
Wiener index of hexagonal systems,  {\it Acta Appl. Math.\/} {\bf
72} (2002) 247--294.

\bibitem{GBS} A. Ghosh, S. Boyd, A. Saberi, Minimizing effective
resitance of graph, {\it SIAM Rev.\/} {\bf 50} (2008) 37--66.

\bibitem{Kl} D.J. Klein, Resistance-distance sum rules, {\it Croat. Chem. Acta\/} {\bf 75} (2002) 633-649.

\bibitem{KR}
D.J. Klein, M. Randi\'{c}, Resistance distance, {\it J. Math. Chem.\/} {\bf 12}
(1993) 81--95.

\bibitem{LNT}
I. Lukovits, S. Nikoli\'{c}, N. Trinajsti\'{c}, Resistance distance
in regular graphs, {\it Int. J. Quantum Chem.\/} {\bf 71} (1999) 217--225.


\bibitem{Pal}
J.L. Palacios, Resistance distance in graphs and random walks, {\it Int.
 J. Quantum Chem.\/} {\bf 81} (2001) 29--33.



\bibitem{YJ} Y. Yang, X. Jiang,  Unicyclic graphs with
extremal Kirchhoff index, {\it MATCH Commun. Math. Comput. Chem.\/}
{\bf 60} (2008) 107--120.

\bibitem{ZD} W. Zhang, H. Deng, The second maximal and minimal Kirchhoff
indices of unicyclic graphs, {\it MATCH Commun. Math. Comput.
Chem.\/} {\bf 61} (2009)  683--695.





\bibitem{ZT1}
B. Zhou, N. Trinajstic, A note on Kirchhoff index, {\it Chem. Phys. Lett.\/} {\bf 445} (2008) 120--123.

\bibitem{ZT2}
B. Zhou, N. Trinajstic, On resistance-distance and Kirchhoff index, {\it J. Math. Chem.\/} {\bf 46} (2009) 283--289.

\bibitem{ZT3}
B. Zhou, N. Trinajstic, Mathematical properties of molecular descriptors based on distances, {\it Croat. Chem. Acta\/} {\bf 83} (2010) 227--242.



\bibitem{ZT4}   B. Zhou, N. Trinajstic, The Kirchhoff index and the matching number, {\it Int. J. Quantum Chem.\/} {\bf 109} (2009) 2978--2981.





%


%





\end{thebibliography}
\end{document}